\documentclass[a4paper,journal,onecolumn,final]{IEEEtran}
\usepackage{showkeys}

\usepackage[noadjust]{cite}

\ifCLASSINFOpdf
  \usepackage[pdftex]{graphicx}
  \graphicspath{{../pdf/}{../jpeg/}}
  \DeclareGraphicsExtensions{.pdf,.jpeg,.png}
\else
  \usepackage[dvips]{graphicx}
	\usepackage{epsfig}
\fi

\usepackage[cmex10]{amsmath}
\usepackage{amssymb}

\usepackage{amsthm}

\usepackage[font=footnotesize]{subfig}
\newtheorem{theorem}{Theorem}
\newtheorem{proposition}[theorem]{Proposition}
\newtheorem{lemma}[theorem]{Lemma}
\newtheorem{corollary}[theorem]{Corollary}
\newtheorem{algorithm}[theorem]{Algorithm}

\theoremstyle{definition}
\newtheorem{definition}[theorem]{Definition}
\newtheorem{example}[theorem]{Example}

\theoremstyle{remark}
\newtheorem{remark}[theorem]{Remark}

\newcommand{\beqa}{\begin{eqnarray*}}
\newcommand{\eeqa}{\end{eqnarray*}}

\newcommand{\field}[1]{\mathbb{#1}}
\newcommand{\bR}{\field{R}}        
\newcommand{\bN}{\field{N}}        
\newcommand{\bC}{\field{C}}        
 \def\cG{\mathcal{G}}

\def\<{\left<}
\def\>{\right>}

\def\inv{^{-1}}

\def\mv1{M_v^1}

\newcommand{\abs}[1]{\lvert#1\rvert}
\newcommand{\norm}[1]{\lVert#1\rVert}

\newcommand{\gab}{\cG (g,\alpha ,\beta )}
\newcommand{\C}{\mathbb{C}}

\newcommand{\N}{\mathbb{N}}
\newcommand{\R}{\mathbb{R}}
\newcommand{\Z}{\mathbb{Z}}

\newcommand{\sprod}[2]{\langle #1, #2 \rangle}

\newcommand{\zweinorm}[1]{\left\|#1\right\|_2}

\newcommand{\bn}{{\boldsymbol{n}}}
\newcommand{\br}{{\boldsymbol{r}}}

\begin{document}
%
\title{Implementation of discretized Gabor frames and their duals}
%
%
%

\author{ Tobias Kloos, Joachim St\"ockler, and Karlheinz Gr\"ochenig
\thanks{This work was supported by the Vienna Science and Technology
  Fund (WWTF) by Grant 10-066 and by the Austrian Science Fund (FWF)
  by Grant P26273-N25.}%
\thanks{K.\ G.\ is with the Faculty of Mathematics,
  University of Vienna, A-1090 Vienna, Austria (email:
  karlheinz.groechenig@univie.ac.at).}%
\thanks{T.\ K.\ and J.\ S.\ are with the Faculty of Mathematics, TU
  Dortmund, D-44221 Dortmund, Germany
  (email: tobias.kloos@tu-dortmund.de; joachim.stoeckler@math.tu-dortmund.de).}%
}
\IEEEpubid{0000--0000/00\$00.00~\copyright~2013 IEEE}


\maketitle

\begin{abstract}
 The usefulness of Gabor frames depends on the easy  computability of
 a suitable dual window. This question is addressed under several
 aspects:  several versions of  Schulz's iterative algorithm for the
 approximation  of the
 canonical dual window are analyzed for their  numerical stability. For
 Gabor frames with totally positive windows or with exponential
 B-splines a direct algorithm yields a family of exact dual windows
 with compact support. It is shown that these dual windows converge exponentially fast
 to the canonical dual window.
\end{abstract}


\section*{Introduction}



The discrete Gabor transform is a useful tool for the analysis and synthesis of nonstationary signals. It is based on the representation of the energy distribution of a signal in the time-frequency plane. Its applications range over the decomposition of musical and acoustical signals \cite{Doerf:2001}, \cite{Russetal:1998}, wireless communication \cite{HlawatschMatz:2011}, \cite{Srirametal:2013} and to the analysis of EEG signals \cite{Blancoetal:1996}, \cite{Chenetal:2010}.
For a given window function $g\in L^2(\R)$ and lattice parameters $\alpha,\beta>0$, the system of all corresponding time-frequency shifts
$$\gab =\{ M_{l\beta}\,T_{k\alpha}\,g = e^{2\pi il\beta\cdot}\,g(\cdot -k\alpha) \mid k,l\in\Z \}$$
is called a Gabor system for $L^2(\R)$. It is called a Gabor frame for
$L^2(\R)$, if  there exist constants $A,B>0$, such that 
\begin{equation}\label{eq:framecond}
   A\|f\|^2 \le \sum_{k,l\in\Z} |\langle f, M_{l\beta} T_{k\alpha} g\rangle|^2
	\le B\|f\|^2,\quad \forall f\in L^2(\R) .
\end{equation}
 The constants $A,B$ are
 called lower and upper frame bounds of $\mathcal{G}(g, \alpha , \beta
 )$. If $\mathcal{G}(g, \alpha , \beta
 )$  fulfills only  the
right hand inequality,  it is called a Bessel sequence and $B$ a Bessel
bound. It is known that the frame inequality~\eqref{eq:framecond}
implies the  existence of  a dual
Gabor frame $\mathcal G(\gamma,\alpha,\beta)$  with dual window
$\gamma\in L^2(\R)$, such that every $f\in L^2(\R)$ can be represented
as
\begin{align}\label{dual:recon}
f = \sum_{k,l\in\Z}\sprod{f}{M_{l\beta}T_{k\alpha}g}\,M_{l\beta}T_{k\alpha}\gamma.
\end{align}
For a given Gabor system, the Gabor transform of a signal $f$ is defined as the analysis operator
\begin{align*}
&\mathcal C_g:L^2(\R)\rightarrow\ell^2(\Z^2),\\
&\mathcal C_gf = (\sprod{f}{M_{l\beta}T_{k\alpha}g})_{k,l\in\Z}.
\end{align*}
The  coefficient $\sprod{f}{M_{l\beta}T_{k\alpha}g}$  represents the
energy distribution of $f$ near the point $(k \alpha , l \beta )$ in
the time-frequency plane. It may also be interpreted as the amplitude
of the frequency $l\beta $ at time $k\alpha $, insofar as such an
interpretation is compatible with the uncertainty principle.
The associated synthesis operator for the reconstructions
\eqref{dual:recon} is the adjoint operator $\mathcal C_g^*$,  and the
frame operator $\mathcal S_g:L^2(\R)\rightarrow L^2(\R)$ is defined by
\begin{align*}
\mathcal S_gf = \mathcal C_g^*\,\mathcal C_gf = \sum_{k,l\in\Z}\sprod{f}{M_{l\beta}T_{k\alpha}g}\,M_{l\beta}T_{k\alpha}g.
\end{align*}

In general, there exist many dual windows suitable for the
reconstruction~\eqref{dual:recon}. The standard choice is the
canonical dual window $\gamma ^\circ = \mathcal{S}_g \inv g$. For a
characterization of all dual windows
see~\cite{Chris:2003,Groech:2001}. Since the applicability and
usefulness of Gabor frames depends heavily on the knowledge and
computability of a dual window,  the numerical construction of dual
windows has motivated numerous studies. As representative
contributions we mention ~\cite{JansSon:2007,strohmer98} and the large
time-frequency analysis toolbox (LTFAT)~\cite{ltfatnote030}.

Our contribution to the analysis of dual Gabor windows is twofold. On a general level, we study  numerically
stable methods for the computation of the canonical dual  window $\gamma^\circ =
\mathcal S_g^{-1}g$. On a specific level,
we study the efficient construction and the behavior of a sequence of
dual windows for Gabor frames with  totally positive window
functions and with exponential B-splines.

We first  present two stable implementations of
a conventional iterative algorithm to approximate $\gamma^\circ$. The
algorithm was originally proposed by Schulz \cite{Schulz:1933} for
matrix inversion. It is based on the Neumann series for the inverse
frame operator and converges quadratically, see Algorithm
\ref{Schulz_iter}. We provide a detailed analysis of the numerical
error and show that our two implementations (the operator version and
the first vector version) are stable. By contrast, the
implementation proposed by Janssen (see \cite{Janssen:2002} and
\cite{Janssen:2003.3}) is often unstable  because  the numerical
error roughly doubles in each step. Therefore  the first two
implementations are much preferable. As two illustrative examples, we
use the window functions MONSTER defined in \cite{JansSon:2007} and a
Gaussian window $g$, in order to compare all three
implementations. The numerical results agree precisely with the
predicted behavior of the numerical error.
\vspace{.3cm}

We then describe recent results on special Gabor
systems whose window function is a totally positive (TP) function of
finite type or an exponential B-spline (EB-spline). TP functions  are
remarkable because so far they are the only window functions for which
a complete characterization of  all lattice parameters such that
$\mathcal{G}(g,\alpha,\beta )$ is a frame is known. More precisely,
 the Gabor system $\gab$, with  a TP
function $g$ of finite type $N\ge 2$ is  a  frame if and only if 
$\alpha\beta<1$~\cite{GroeStoe:2013}. Subsequently, similar arguments
in~\cite{KloStoe:2014} showed   that the Gabor
system $\mathcal G(B_\Lambda,\alpha,\beta)$ of an EB-spline
constitutes a frame for $\alpha = 1$, $\beta<1$,  and some other lattice
parameters, too. The proofs also provide a constructive method for the
computation of infinitely many dual windows $\gamma_L$ with compact
support, which we summarize in Algorithm~\ref{algo:gamma} for TP
functions of finite type and Algorithm~\ref{algo2:gamma} for
EB-splines.

This construction offers several new and useful aspects that are special
for TP windows and not shared by general window functions.

(i) Algorithms~\ref{algo:gamma} and ~\ref{algo2:gamma} provide a
family of dual windows $\gamma _L$ both in finite and infinite
dimensional models, namely for
continuous signals in $L^2(\R )$, for discrete signals in $\ell ^2(\Z
)$, and for periodic discrete signals in $\C ^N$. Currently available
toolboxes, such as LTFAT~\cite{ltfatnote030},  work only for finite-dimensional
signals.

(ii) The dual windows $\gamma _L$ possess compact support of size
$\mathcal{O}(L)$, whereas the canonical dual $\gamma ^\circ $ is known
to have infinite support.

(iii) The dual windows $\gamma _L$ are exact and
satisfy~\eqref{dual:recon}. This is  in contrast to the standard
iterative methods for the approximation of the canonical dual (see
e.g. Algorithm~\ref{Schulz_iter}), which generate only approximations
of a dual  window. 

As our main mathematical result
we prove that the dual windows
$\gamma_L$ are good approximations of the canonical dual window
$\gamma^\circ = \mathcal S_g^{-1}g$ and we show that they converge
exponentially fast to the canonical dual window, i.e., 
$\|\gamma_L-\gamma^\circ\|_2=\mathcal{O}(e^{-\rho L})$.
Therefore, by specifying the
parameter $L$, Algorithms~\ref{algo:gamma} and \ref{algo2:gamma}
provide a dual window $\gamma_L$ with compact support and which
approximates the canonical dual at a desired rate. 

 The proof uses some ideas of the non-symmetric finite
section method, but also requires a new technique related to the
formulation of the Moore-Penrose pseudo-inverse of infinite matrices
in terms of orthogonal projections.
\vspace{.3cm}

As our main numerical contribution, we study and implement
the case of discrete Gabor frames. We present some fast and stable
algorithms to evaluate and discretize TP functions and EB-splines and
their  dual windows
computed by the Algorithms~\ref{algo:gamma} and
\ref{algo2:gamma}. These algorithms are proposed as extensions to the
Large Time Frequency Analysis Toolbox described in
\cite{ltfatnote030}.

The paper is organized as follows: In section~\ref{sec:hotelling} we
study the numerical stability of a fast iterative  algorithm for the
approximation of the canonical dual window. In section~\ref{sec:known}
we summarize the algorithms for the construction of dual windows of TP
functions and EB splines. In section~\ref{sec:approx} we formulate and
discuss the main theorems about the convergence of the compactly
supported dual windows $\gamma _L$ to the canonical dual window $\gamma
^\circ $. Section~\ref{sec:discrete} explains some details about the
implementation of Gabor frames with TP functions and EB splines.
The appendix contains the technical details of the
proofs of the main results.


\section{Some iterative algorithms for approximating the canonical dual}\label{sec:hotelling}

In this section we describe two iterative algorithms for approximating
the canonical dual of an arbitrary frame $\mathcal F=\{f_j\}_{j\in I}$
for a Hilbert space $\mathcal H$. The central part of such algorithms
is the approximation of the inverse of the corresponding frame
operator
$$S_{\mathcal F}:\mathcal H\to \mathcal H,\qquad
S_{\mathcal F} h=\sum_{j\in I} \langle h,f_j\rangle f_j.$$
We  discuss the convergence and the
numerical stability of various implementations. Finally we present
some numerical tests. The following approximation schemes are proposed
in the literature.
\begin{algorithm}[Frame algorithm]
Choose $0<\lambda<2/B$, with $B$ the upper frame bound of $\mathcal F$. Then $q:=\norm{I-\lambda\,S_{\mathcal F}}<1$ and
$$S_{\mathcal F}^{-1} = \lambda\,\sum_{n=0}^\infty (I-\lambda\,S_{\mathcal F})^n.$$
The partial sums of this Neumann series can be computed iteratively  by
\begin{align}
&K_0 = \lambda\,I, \notag \\
&K_{k+1} = \lambda\,I + (I-\lambda\,S_{\mathcal F})\,K_k,\ \
k\in\N_0. \label{eq:cmr}
\end{align}
The convergence rate is  of order $\norm{S_{\mathcal F}^{-1}-K_{k}} =
\mathcal O(q^k)$.
The $k$'th approximation of the canonical dual frame $\mathcal F^\circ
= \{S_{\mathcal F}^{-1}\,f_j\}_{j\in I}$ is given by
$\{K_k\,f_j\}_{j\in I}$.
\end{algorithm}
This algorithm is very robust, but slow if $q$ is close to one.
This algorithm can be accelerated~\cite{gro93ieee} with
conjugate gradient techniques and with a convergence rate
$\big(\frac{\sqrt{B}- \sqrt{A}}{\sqrt{B}+ \sqrt{A}}\big)^k$ after $k$
iterations regardless of whether the frame bounds $A,B$ are known or
not.

An even faster method goes back to Schulz~\cite{Schulz:1933} and
Hotelling~\cite{Hotelling:1942}.

\begin{algorithm}[Schulz iteration]\label{Schulz_iter}
 Choose  $0<\lambda < 2/B$.  The version of Schulz iteration  with
 "initial scaling"~\cite[Algorithm~IV]{Janssen:2003.3} is
\begin{align}\notag
&J_0 = \lambda\,I,\\\label{Schulz_op}
&J_{k+1} = 2\,J_k - J_k\,S_{\mathcal F}\,J_k,\ \ k\in\N_0.
\end{align}
This iteration implies  the identity $J_k = K_{2^k-1}$ and is
therefore connected  to the frame algorithm. The Schulz algorithm  converges quadratically, i.e.
\begin{equation}
  \label{eq:c13}
\norm{S_{\mathcal F}^{-1}-J_{k+1}} \le \norm{S_{\mathcal
    F}}\,\norm{S_{\mathcal F}^{-1}-J_k}^2 = \mathcal{O}(q^{2^{k+1}}) .
\end{equation}
\end{algorithm}
This algorithm was first described by Schulz \cite{Schulz:1933}  who used this method   for matrix inversion.
\begin{IEEEproof}
The claims in Algorithm~\ref{Schulz_iter} are proved by induction. Since Schulz's algorithm is not as known as other iterative
algorithms, we sketch the main steps.
We first show that
\begin{equation}
  \label{eq:c19}
  I - S_\mathcal{F}J_k  = I - J_k S_\mathcal{F} = (I - \lambda S_\mathcal{F}
  )^{2^k} \, .
\end{equation}
Assuming that \eqref{eq:c19} is correct for $k\in \mathbb{N}_0$, we
obtain
\begin{align*}
  I - S_\mathcal{F}J_{k+1} &= I - S_\mathcal{F} (2 J_k - J_k S_\mathcal{F} J_k) \\
&= (I - S_\mathcal{F}J_k)^2 = \Big( (I - \lambda S_\mathcal{F} )^{2^k} \Big)^2
\, ,
\end{align*}
as claimed.
Using~\eqref{eq:c19}, we show again by induction that  $J_k = K_{2^k-1} =\lambda \sum
_{j=0}^{2^k-1} (I - \lambda S_\mathcal{F} )^j $:
\begin{align*}
\lambda & \sum_{j=0}^{2^{k+1}-1}  (I-\lambda S_\mathcal{F} )^j  \\
                    &           = \lambda
                               \sum_{j=0}^{2^{k}-1}(I-\lambda S_\mathcal{F}
                               )^j + (I-\lambda S_\mathcal{F} )^{2^k} \, \lambda
                               \sum_{j=0}^{2^{k}-1}(I-\lambda S_\mathcal{F} )^{j}\\
                               &= J_k +  (I - J_k S_\mathcal{F} ) J_k =
                               J_{k+1} \, .
\end{align*}
The quadratic convergence rate now  follows  from the convergence
properties of the Neumann series.
\end{IEEEproof}

%

We discuss the implementation of Algorithm~\ref{Schulz_iter} in the
case of a Gabor frame $\mathcal G(g,\alpha,\beta)$. Recall that the
canonical dual frame is determined by the dual window
$\gamma^\circ=\mathcal S_g^{-1}g$. We compare three different
implementations of the Schulz iteration  and provide some heuristics
for their  numerical
stability.
\medskip

\noindent (i) \emph{Operator form}:  The numerical computation of the
  Schulz iteration as stated in Algorithm~\ref{Schulz_iter} provides
  operators $\hat{J}_k = J_k + E_k$, where $E_k$ denotes the
  accumulated forward error. Let  $Y_{k+1}$ denote the new roundoff
  error in the $k+1$'st iteration, then the operator  after  $k+1$
  iterations of  (\ref{Schulz_op}) is
\begin{align*}
\hat{J}_{k+1} &= 2\,\hat{J}_k - \hat{J}_k\,\mathcal S_g\,\hat{J}_k + Y_{k+1}\\
              &= J_{k+1} + E_k\,(I-\mathcal S_g\,J_k) +
              (I-J_k\,\mathcal S_g)\,E_k
              + Y_{k+1} +
              \mathcal O(\norm{E_k}^2).
\end{align*}
 Since $I-\mathcal S_g\,J_k = I-J_k\,\mathcal S_g =
(I-\lambda\,\mathcal S_g)^{2^k}$,  we have
$$\norm{E_{k+1}} = \norm{\hat{J}_{k+1} - J_{k+1}}  \le 2\,q^{2^k}\,\norm{E_k} + \norm{Y_{k+1}} + \mathcal O(\norm{E_k}^2).$$
This estimate shows that the error accumulated in the first $k$
iterations is damped and only a new round-off error is added.  Hence this iteration is numerically stable.
\medskip

\noindent (ii) \emph{ Vector form}:  We compute approximations of the
  dual window $\gamma^\circ$ directly by setting $\gamma_k := J_k\,g$
  with the following algorithm:
\begin{align}\notag
&\gamma_0    = \lambda\, g \\\label{Hotelling_el}
&\gamma_{k+1}= 2\, \gamma_k - \mathcal C^*_{\gamma_k}\,\mathcal
C_{g}\,\gamma_k,\ \ k\in\N_0 \, .
\end{align}
\begin{IEEEproof}
We use  the fact that $S_\mathcal{F}$ and thus all
$J_k$ commute with the time-frequency shifts $M_{l\beta } T_{j\alpha
}$. This commutation rule implies  the  identity
\begin{align*}
\mathcal C_{\gamma_k}^*\,\boldsymbol{c} = \sum_{j,l\in
  \Z}c_{j,l}\,J_k\,M_{l\beta}\,T_{j\alpha}\,g = J_k\,\mathcal
C_g^*\,\boldsymbol{c}
\end{align*}
for all $\boldsymbol{c}\in\ell^2(\Z^2)$.  Consequently
 \begin{align*}
   \gamma _{k+1} &= J_{k+1} g = 2 J_k g -  J_k \mathcal{S}_g J_kg \\
   &= 2\gamma _k -  J_k \mathcal{C}^*_g \mathcal{C}_g \gamma_k = 2\gamma _k
   -  \mathcal{C}^*_{\gamma _k} C_g \gamma _k \, .
 \end{align*}
\end{IEEEproof}
 The numerical computation yields  $\hat{\gamma}_k = \gamma_k + e_k$,
 where $e_k$ denotes the accumulated forward error. Let  $y_{k+1}$
 denote the new roundoff error, then in  the $k+1$'th
 step of the iteration (\ref{Hotelling_el})  we have
\begin{align*}
\hat{\gamma}_{k+1} &= 2\,\hat{\gamma}_k - \mathcal C^*_{\hat{\gamma}_k}\,\mathcal C_{g}\,\hat{\gamma}_k + y_{k+1}\\
              &= \gamma_{k+1} + (e_k-\mathcal C_{\gamma_k}^*\,\mathcal
              C_g \,e_k) + (e_k - \mathcal  C_{e_k}^*\,\mathcal
              C_g\,\gamma_k)
              + y_{k+1} + \mathcal
              O(\norm{e_k}_2^2) \, .
\end{align*}
 The aforementioned estimates give
$$\norm{(I-\mathcal C_{\gamma_k}^*\,\mathcal C_g )\,e_k}_2 = \norm{(I-J_k\,\mathcal S_g)\,e_k}_2\le q^{2^k}\,\norm{e_k}_2.$$
Moreover the Janssen representation (see \cite[p.131]{Groech:2001}) gives
$$\mathcal C_{e_k}^*\,\mathcal C_g\,\gamma_k = \frac{1}{\alpha\beta}\sum_{j,l\in\Z}\sprod{e_k}{M_{j/\alpha}\,T_{l/\beta}\,g}M_{j/\alpha}\,T_{l/\beta}\,\gamma_k.$$
The last expression, when $\gamma_k$ is replaced by $\gamma^\circ$, is
the orthogonal projection $\Pi_{V_g} e_k$ of $e_k$ onto $V_g :=
\overline{\mathrm{span}}\bigl(\mathcal G(g,1/\beta,1/\alpha)\bigr)$
and therefore
$$\norm{e_k - \mathcal  C_{e_k}^*\,\mathcal
  C_g\,\gamma^\circ}_2\le\norm{e_k}_2$$
and
\begin{align*}
\norm{e_k &- \mathcal  C_{e_k}^*\,\mathcal
  C_g\,\gamma _k}_2 \le \norm{e_k - \mathcal  C_{e_k}^*\,\mathcal
  C_g\,\gamma^\circ}_2
  +\norm{C_{e_k}^*\,\mathcal
  C_g\,(\gamma^\circ - \gamma _k)}_2 \leq \norm{e_k}_2 +
\mathcal{O}(q^{2^k}) \, .
\end{align*}
Since $\norm{\gamma^\circ - \gamma_k}_2=\mathcal O(q^{2^k})$ by
\eqref{eq:c13},  the new numerical error is
$$\norm{e_{k+1}}_2 \leq (1+q^{2^k}) \norm{e_k}_2 +
  \norm{y_{k+1}}_2 + \mathcal{O}(\norm{e_k}^2_2 + q^{2^k})  \, .
$$
In contrast to the operator version, $\norm{e_k}_2$ enters linearly
with coefficient $\approx 1$. Thus
the numerical stability is plausible and  is also
confirmed by our numerical results in Example~\ref{numerical_results}.
\medskip

\noindent (iii) \emph{Janssen's alternative vector version}:  Janssen \cite[Algorithm~IV]{Janssen:2003.3}
  proposes  the
  approximations
\begin{align}\notag
&\gamma_0    = \lambda\, g, \\\label{Janssen_el}
&\gamma_{k+1}= 2\, \gamma_k - \mathcal C^*_{\gamma_k}\,\mathcal C_{\gamma_k}\,g,\ \ k\in\N_0.
\end{align}
However,  in their numerical tests of \eqref{Janssen_el}   Janssen and S{\o}ndergaard \cite{JansSon:2007}   observed some numerical instability. With notations as in (ii), the numerical error in step $k+1$ is
\begin{align*}
  e_{k+1}&=\hat{\gamma}_{k+1}-\gamma_{k+1} \\
&=2\,e_k - \mathcal
  C_{\gamma_k}^*\,\mathcal C_{e_k}\,g - \mathcal C_{e_k}^*\,\mathcal
  C_{\gamma_k}\,g  + y_{k+1} + \mathcal{O}(\norm{e_k}_2^2).
\end{align*}
We show that the error may grow by at least a factor of $2$ in each step. Note that $V_g$ in (ii) is a proper subset of $L^2(\R)$, if $\alpha\beta<1$. Hence $\norm{(I-\Pi_{V_g})\,e_1}_2\approx \varepsilon$, where $\varepsilon$ denotes the machine accuracy. The Janssen representation implies, that $\mathcal C_{\gamma_k}^*\,\mathcal C_{e_k}\,g$ and $\mathcal C_{e_k}^*\,\mathcal C_{\gamma_k}\,g$ are in $V_g$. Therefore,
\begin{align*}
(I-\Pi_{V_g})\,e_{k+1} = 2\,(I-\Pi_{V_g})\,e_k&+(I-\Pi_{V_g})\,y_{k+1}
   +\mathcal O(\norm{e_k}_2^2)
\end{align*}
This shows, that the error components in $V_g^\perp$ can double in
each step.  In Example~\ref{numerical_results} we demonstrate that
this numerical instability  may indeed occur.

\begin{center}
\begin{figure}
\hspace*{4.5cm}%
\includegraphics[width=8.8cm]{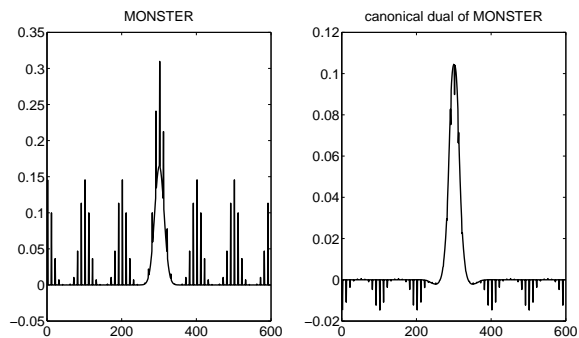}
\caption{MONSTER function and its canonical dual for $\alpha=20$ and $\beta=1/50$, computed by a routine of LTFAT \cite{ltfatnote030}.}\label{plot:monster}
\end{figure}
\end{center}

\begin{center}
\begin{figure}
\hspace*{4.5cm}%
\includegraphics[width=8.8cm]{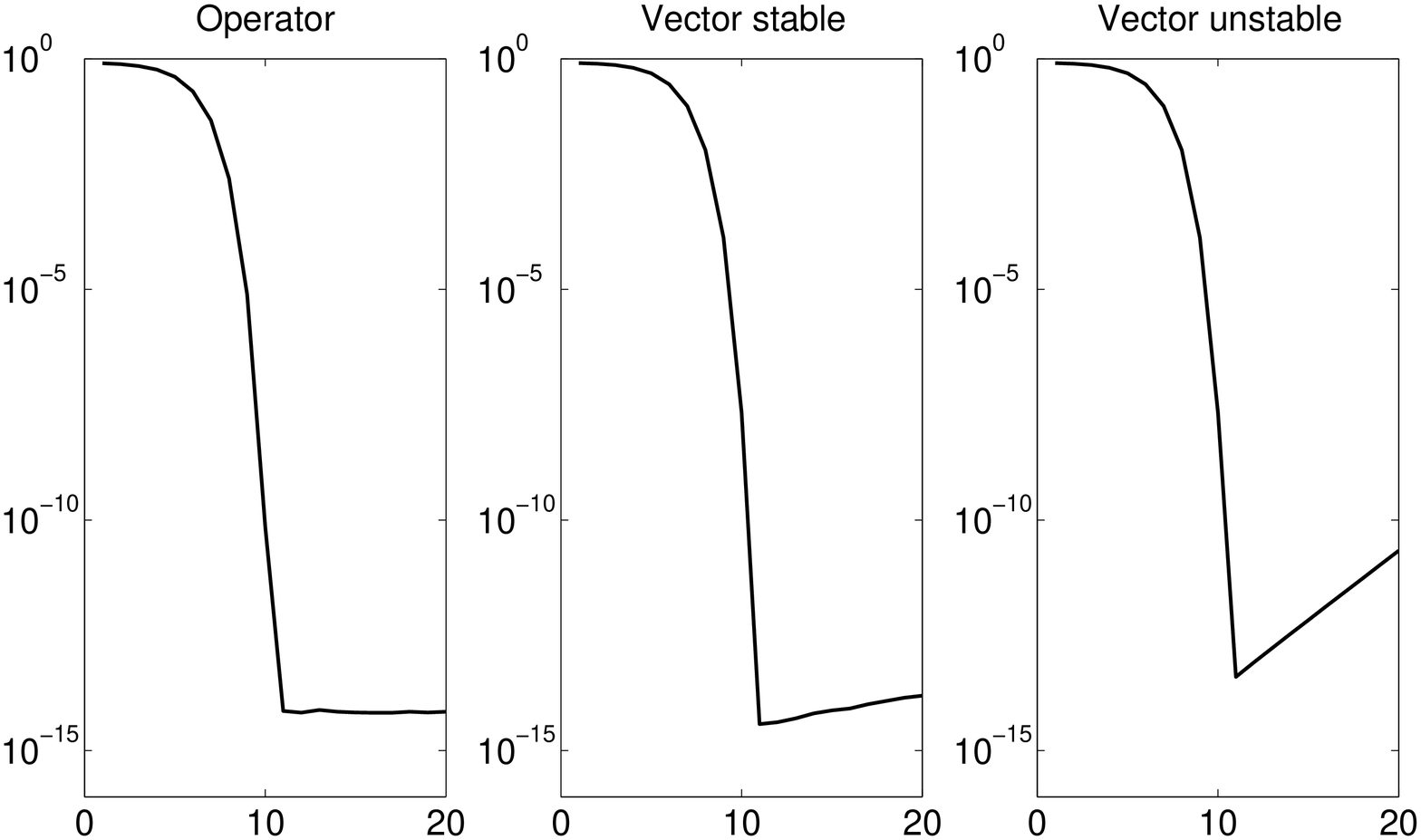}

\caption{Error $\norm{\mathcal{C}_{\gamma^\circ}-\mathcal{C}_{\hat\gamma_k}}$ of the three implementations of Algorithm \ref{Schulz_iter} for $20$ iteration steps to approximate the canonical dual window of MONSTER for $\alpha=20$ and $\beta=1/50$.}\label{plot:error}
\end{figure}
\end{center}

\begin{center}
\begin{figure}[h]
\hspace*{4.5cm}%
\includegraphics[width=8.8cm]{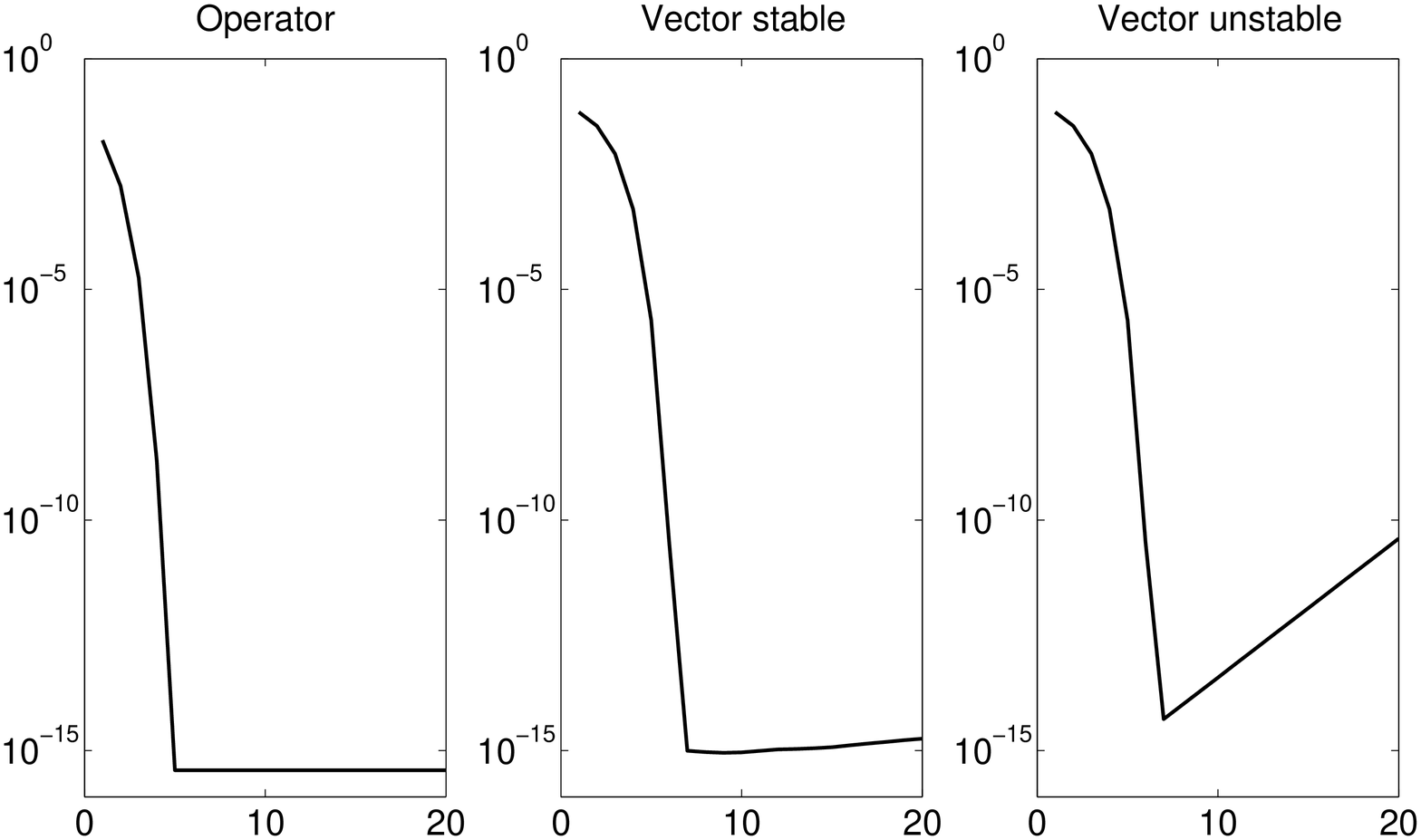}

\caption{$\norm{\mathcal{C}_{\gamma^\circ}-\mathcal{C}_{\hat\gamma_k}}$ of the three implementations of Algorithm \ref{Schulz_iter} for $20$ iteration steps to approximate the canonical dual window of the Gaussian $g(x) = e^{-\pi x^2/600}$ for $\alpha=20$ and $\beta=1/50$.}\label{plot:error2}
\end{figure}
\end{center}

\begin{example}\label{numerical_results}
As it was described in (iii) the implementation (\ref{Janssen_el}) by
Janssen can have some stability problems and should be applied
carefully. We use the function MONSTER (see Figure~\ref{plot:monster})
in \cite{JansSon:2007} and $\alpha=20$, $\beta=1/50$ for our numerical
tests of all three implementations of the Schulz iteration. As we can
see in Figure~\ref{plot:error}, the operator version is stable, the
vector version of (ii) is also useful, while the error of the
implementation in (iii) explodes. The same conclusions hold for the window function  $g(x) =
e^{-\pi x^2/600}$ with  $\alpha= 20, \beta= 1/50$,  as is shown in Figure~\ref{plot:error2}.
\end{example}

\newpage
\section{Gabor frames of totally positive functions and exponential B-splines}
\label{sec:known}

 We now consider totally positive functions of finite type and
 exponential B-splines as window functions. The TP functions attracted
 much interest recently, as they provide new examples of window
 functions for which the necessary density condition $\alpha\beta<1$
 of the lattice parameters is also sufficient
 \cite{GroeStoe:2013}. For detailed information on total positivity of
 functions and matrices see \cite{Karl:1968}, for a detailed
 introduction to exponential B-splines see \cite{Schum:1981} and
 \cite{deVLor:1993}.

\begin{definition}{\bf \cite{Schoenberg:1947}, \cite{Schoenberg:1951}}
An integrable function $g:\R\rightarrow\R$ is called totally positive
(TP), if
its Fourier transform factors as
\begin{align}\notag
\hat{g}(\omega) &=  \int_{-\infty}^{\infty} g(t) e^{-2\pi it\omega}\,dt \\ \label{TPall}
&= C e^{-\eta\omega^2} e^{-2\pi i\delta\omega}
\prod_{\nu=1}^{\infty}\frac{e^{2\pi i\delta_{\nu}\omega}}{1+2\pi i\delta_{\nu}\omega},
\end{align}
where $C,\eta,\delta,\delta_{\nu}$ are real parameters with
$$C>0,\ \ \eta\geq0, \
0<\eta+\sum_{\nu=1}^{\infty}\delta_{\nu}^2<\infty\, .$$
\end{definition}

Note that in Schoenberg's terminology there also exist TP functions
that  are not integrable. Therefore the given definition  does not
include the most general case of TP functions. Here we only consider TP
functions $g_\Delta$ of finite type $N\in\N$, with $\eta =0$, $\Delta=(\delta_1,\ldots,\delta_N)$ and
\begin{equation}\label{eq:TPFFT}
    \widehat{g_\Delta}(\omega) = \prod_{\nu = 1}^N(1+2 \pi i\delta_\nu\omega)\inv.
\end{equation}
For $N=1$ and $\Delta=(\delta)$ we have the one-sided exponential
$$
  g_\delta(x)=\frac{1}{|\delta|} \,
  e^{-\delta\inv x}\, \chi_{(0,\infty)}(\delta x)\quad\hbox{for}\quad x\in\R\setminus\{0\},
$$
with support $[0,\infty)$ for positive $\delta$ and $(-\infty,0]$ for negative $\delta$, and for $N>1$ and $\Delta=(\delta_1,\ldots,\delta_N)$, $g_\Delta$ is the $N$-fold convolution
$$
   g_\Delta=
   g_{\delta_1}\ast\cdots \ast
   g_{\delta_N}.
$$
Especially, if $\delta_\nu=\delta\ne 0$ for all $1\le \nu\le N$, then
$$
  g_\Delta(x)= \frac{|x|^{N-1}}{|\delta|^N(N-1)!}\,
  e^{-\delta\inv x}\, \chi_{(0,\infty)}(\delta x).
$$
The functions $g_\Delta$  are nonnegative, have infinite support,  and exponential
decay, precisely, if $\tau = (\max\{\,2\pi\,\abs{\delta_\nu} \mid \nu
= 1,\ldots,N\})^{-1}$ and $\epsilon >0$, then there is a constant $c =
c(\epsilon )$ such that
\begin{align}\label{TP:decay}
   0\le g_\Delta (x)\le c\,e^{-(\tau-\varepsilon)\abs{x}}.
\end{align}
Moreover, if $\delta_1\ne \delta_N$ and
we specify $g_{\delta}(0)=1/|2\delta|$,
the recurrence relation
$$g_{\delta_1,\ldots,\delta_N}=\frac{\delta_1\inv g_{\delta_2,\ldots,\delta_N}-
\delta_N\inv g_{\delta_1,\ldots,\delta_{N-1}}}{\delta_1\inv-\delta_N\inv}
$$
holds.
For later use, we denote by $m$ (resp. $n$) the number of positive (resp. negative) parameters $\delta_\nu$.

The main result in \cite{GroeStoe:2013} shows that the Gabor system $\gab$ of a TP function of finite type $m+n\ge 2$ constitutes a Gabor frame if and only if $\alpha\beta<1$. The proof provides an algorithm for the computation of a dual window $\gamma$ with compact support. This method was adapted in \cite{BanGroeStoe:2013}, \cite{Klo:2012} in order to supply infinitely many dual windows $\gamma_L$. The computation of $\gamma_L(x+j\alpha)$, with $j\in\Z$, is performed by computing a single row of a left-inverse of the biinfinite pre-Gramian matrix
\begin{align}\label{pregram}
P_g(x) := \left( \overline{g}\left(x+\alpha j -\tfrac{k}{\beta}\right) \right)_{j,k\in\Z}.
\end{align}
The structure of the left-inverse of $P_g(x)$ heavily depends on the property that $g$ is totally positive. We include the algorithm for the reader's convenience.

\begin{algorithm} {\bf \cite{BanGroeStoe:2013}}
\label{algo:gamma}
Input parameters are the parameter vector $\Delta = (\delta _1,
\delta _2, \dots , \delta _N)$ of the
window $g$, the lattice parameters $\alpha ,\beta > 0$ with
$\alpha\beta<1$,
a parameter $L\in\bN_0$ controlling the support size of the dual window $\gamma_L$,
and a point $x\in[0,\alpha)$.

Output parameters are integers $i_1(L),i_2(L)$ and
the vector of values $\gamma_L(x+\alpha j)$,
$i_1(L)\le j\le i_2(L)$, in the support of $\gamma_L$, such that $\gamma_L:\R\rightarrow\R$ defines a dual window of $g$.

  \begin{enumerate}
  \item Set $r:=\left\lfloor \frac{1}{1-\alpha \beta} \right\rfloor$, $ k_1 = -(r+1)m$ and
      $k_2 = (r+1)n$.
  \item Set $ k_1(L) := k_1-L$ and $k_2(L) := k_2+L$,
    \begin{align*}
    i_1(L) &:= \left\lfloor \frac{k_1-L + m -
        1}{\alpha\beta}-\frac{x}{\alpha}\right\rfloor+1,\\
    i_2(L) &:=
    \left\lceil \frac{k_2+L - n +
        1}{\alpha\beta}-\frac{x}{\alpha}\right\rceil-1.
    \end{align*}
   \item Set  $ P_L(x) := (p_{j,k})_{i_1(L)\le j \le i_2(L),~ k_1(L)\le k \le k_2(L)}$, where
    \[  p_{j,k}= g \left(x+\alpha j- \frac{k}{\beta}\right).
    \]
  \item Compute the
	pseudoinverse $$P_L(x)^\dagger=(q_{k,j})_{k_1(L)\le k \le k_2(L),~ i_1(L)\le j \le i_2(L)}$$
    of $P_L(x)$.
  \item Take the row with index \(k=0\) of $P_L(x)^\dagger$. Its
    coefficients define the values of the dual window $\gamma_L$ at the
    points $\{x+\alpha j \,|\,i_1(L) \le j \le i_2(L)\}$, i.e.
    \begin{equation}\label{eq:set-gamma}
    \gamma_L(x+\alpha j):=\begin{cases} \beta\, q_{0,j}&\text{, if } i_1(L)\le j \le i_2(L),\\
		0&\text{, if } j<i_1(L)~\text{or } j>i_2(L).
		\end{cases}
    \end{equation}
  \end{enumerate}
\end{algorithm}
In particular, the support of the dual window $\gamma _L$ is contained
in the interval $[\alpha i_1(L), \alpha (i_2(L)+1)]$ of length of the
order $ \beta\inv (k_2 - k_1 +2L )$. Thus the parameter $L$
labels the size of the support.

A  related class of window functions  is the class of exponential B-splines. These functions are positive and have compact support, a property which is desirable in some applications.

\begin{definition}
For $\Lambda=(\lambda_1,\ldots,\lambda_N)\in\R^N$ the exponential B-spline (EB-spline) $B_\Lambda$ with knots $0,1,\ldots,N$ is given by its Fourier transform
\begin{align}\label{eq:EBspline}
\widehat{B_\Lambda}(\omega) = \prod_{\nu=1}^N \frac{e^{\lambda_\nu-2\pi i\omega}-1}{\lambda_\nu-2\pi i\omega}.
\end{align}
\end{definition}

 In \cite{KloStoe:2014}, it was shown that the Gabor system $\mathcal
 G(B_\Lambda,\alpha,\beta)$ of every EB-spline constitutes a frame for
 $\alpha = 1$, $\beta<1$ (and also  some other lattice
 parameters). Similar to the case of TP functions, the following
 algorithm provides dual windows $\gamma_L$ of the Gabor frame
 $\mathcal{G}(B_\Lambda,\alpha,\beta)$.

\begin{algorithm}~{\bf \cite{KloStoe:2014}}
\label{algo2:gamma}
Input parameters are the parameter vector $\Lambda\in\R^N$ of the
window $B_\Lambda$, the lattice parameters $\alpha ,\beta > 0$,
a parameter $L\in\bN_0$ controlling the support size of the dual window $\gamma_L$,
and a point $x\in\left[\frac{N-\alpha}{2},\frac{N+\alpha}{2}\right)$.

Output parameters are integers $i_1(L),i_2(L)$ and
the vector of values $\gamma_L(x+\alpha j)$,
$i_1(L)\le j\le i_2(L)$, in the support of $\gamma_L$, such that $\gamma_L:\R\rightarrow\R$ defines a dual window of $B_\Lambda$.

  \begin{enumerate}
  \item Set $ k_2(L) = \left\lfloor \frac{N\beta+\alpha\beta}{2(1-\alpha\beta)} \right\rfloor +1 +L$ and $k_1(L) = -k_2(L)$.
  \item Set $i_1(L)\le k_1(L)$ and $i_2(L)\ge k_2(L)$, such that
	\begin{align*}
	B_\Lambda\left(x+(i_1(L)-1)\alpha-\tfrac{k_1(L)-1}{\beta}\right)\neq 0,\\
	B_\Lambda\left(x+i_1(L)\alpha-\tfrac{k_1(L)-1}{\beta}\right) = 0,\\
	B_\Lambda\left(x+(i_2(L)+1)\alpha-\tfrac{k_2(L)+1}{\beta}\right)\neq 0,\\
	B_\Lambda\left(x+i_2(L)\alpha-\tfrac{k_2(L)+1}{\beta}\right) = 0.
	\end{align*}
	\item Set  $ P_L(x) := (p_{j,k})_{i_1(L)\le j \le i_2(L),~ k_1(L)\le k \le k_2(L)}$, where
    \[  p_{j,k}= B_\Lambda \left(x+\alpha j- \frac{k}{\beta}\right).
    \]
  \item Compute the
	pseudoinverse $$P_L(x)^\dagger=(q_{k,j})_{k_1(L)\le k \le k_2(L),~ i_1(L)\le j \le i_2(L)}$$
    of $P_L(x)$.
  \item Take the row with index \(k=0\) of $P_L(x)^\dagger$. Its
    coefficients define the values of the dual window $\gamma_L$ at the
    points $\{x+\alpha j \,|\,i_1(L) \le j \le i_2(L)\}$, i.e.
    \begin{equation}\label{eq:set-gamma2}
    \gamma_L(x+\alpha j):=\begin{cases} \beta\, q_{0,j}&\text{, if } i_1(L)\le j \le i_2(L),\\
		0&\text{, if } j<i_1(L)~\text{or } j>i_2(L).
		\end{cases}
    \end{equation}
  \end{enumerate}
\end{algorithm}

We would like to emphasize the following point: (i)  These algorithms determine the precise values of a
dual window $\gamma $ with compact support in $L^2(\bR )$ and not just
a discrete approximation in a finite-dimensional vector space, as is
done in most existing  algorithms
in~\cite{ltfatnote030,strohmer98}.

(ii) Although the problem of finding a dual window is by nature
infinite-dimensional in $L^2(\R )$, the computation of $\gamma _L$  on
a grid $\tfrac{\alpha}{M} \Z $ requires only the pseudo-inversion of
$M$ finite-dimensional matrices.

\vspace{.5cm}


\section{Approximation of the canonical dual of TP functions and EB-splines}\label{sec:approx}

We fix the parameters $\alpha,\beta>0$ of the Gabor frame and let $g$
be either a TP function of finite type with parameter set $\Delta$ or
an EB-spline with parameter set $\Lambda$. The pre-Gramian matrix
\eqref{pregram} plays a central role in the characterization of the
frame bounds $A,B$ of the Gabor system $\gab$, and in finding dual
windows $\gamma$.

In this section we address the question how the sequence of dual
windows $\gamma _L$ of Algorithms~\ref{algo:gamma} and
~\ref{algo2:gamma} are related to the canonical dual window $\gamma
^\circ$.
Our main goal  is to prove that the compactly supported dual windows
$\gamma_L$ computed by the  Algorithms~\ref{algo:gamma} and
\ref{algo2:gamma} from  \cite{BanGroeStoe:2013} and
\cite{KloStoe:2014}  approximate the canonical dual window $\gamma^\circ$ at a rate
\[
   \zweinorm{\gamma_L- \gamma^\circ} \le \tilde{c}\,e^{-\rho L},
\]
where $L$ determines the support of $\gamma_L$ and $\rho>0$.

As a first step we show that the norm of $\gamma _L$ remains bounded
as $L $ tends to $\infty $. 

\begin{theorem}
  \label{thm:bounds-PL}
Let $g$ be a TP function of finite type as defined by
(\ref{eq:TPFFT}). Then there exist constants $A,B>0$ independent of
$L\in \N _0$ such that
\begin{align}\label{PLbounds}
   A \|\boldsymbol{c}\|_2 \le \|  P_L(x) \boldsymbol{c}\|_2 \le    B \|\boldsymbol{c}\|_2
\end{align}
for all $\boldsymbol{c}\in \bC^{k_2(L)-k_1(L)+1}$ and $x\in [0,\alpha)$, where $P_L(x)$ is the finite section of the corresponding biinfinite pre-Gramian matrix $P_g(x)$ as described in Algorithm~\ref{algo:gamma}. Consequently, $\zweinorm{\gamma_L}\leq \sqrt{\alpha}\,\beta\,A^{-1}$, for all $L\in\N$.
\end{theorem}

The precise proof will be given in the appendix. The proof idea goes
as follows:
The upper bound $B$ is easily  obtained from Schur's test based on
the exponential decay of $g$ in (\ref{TP:decay}). For the lower bound
$A$, we choose $2L+1$ submatrices $R_\ell$ of $P_L(x)$, with
$k_2-k_1+1$ columns each, and build a left-inverse $Q_L(x)$ of
$P_L(x)$ by the selection of specific rows of $R_\ell^\dagger$. The
result in \cite[Theorem~9]{GroeStoe:2013} shows that all matrices
$R_\ell^\dagger$ are bounded uniformly in $x$. Their upper bound $C$
will provide the upper bound
$A^{-1}=C\,\tfrac{k_2-k_1+1}{\alpha\beta}$ for all matrices $Q_L(x)$
uniformly in $x$ and $L$. Hence $A$ is a suitable constant for the
lower bound of $P_L(x)$. 

To prove the rate of approximation of the dual windows $\gamma_L$, we
recall   that the canonical dual window $\gamma^\circ = \mathcal
S_g^{-1}g$ can be expressed by the Moore-Penrose pseudoinverse of
$P_g(x)$, namely
\begin{equation}
  \label{eq:c15}
\gamma^\circ(x+\alpha j) = \beta\,(P_g(x)^\dagger)_{0,j}   
\end{equation}
for all $j\in\Z$ and all $x\in[0,\alpha)$. This is a direct consequence of the
Wexler-Raz criterion for the dual windows of $g$ \cite[Theorem~7.3.1]{Groech:2001}, and the minimal $L^2$-norm of
the canonical dual among all duals of $g$ \cite{Janssen1995}.
We will therefore show that the zeroth row of $P_L(x)^\dagger$ approximates the zeroth
row of $P_g(x)^\dagger$ at an exponential rate. For this we need  a
new result on non-symmetric finite sections of biinfinite
matrices. The following theorem is  not covered by the results in
\cite{GroeRzeStroh:2010} and may be of independent interest.  To fix the notation, for $\boldsymbol{n} =
(n_1,n_2) \in\N^2$ and $b\in\ell^2(\Z)$, we let $\mathcal
P_{\boldsymbol{n}}$ with
$$\mathcal P_{\boldsymbol{n}}\,b = (\ldots,0,b_{-n_1},b_{-n_1+1},\ldots,b_{n_2-1},b_{n_2},0,\ldots)^T$$
be the orthogonal projection onto the $n_1+n_2+1$-dimensional subspace
$\mathcal P_{\boldsymbol{n}}\,\ell^2(\Z) \cong \C^{n_1+n_2+1}$. For a
biinfinite matrix $U=(u_{j,k})_{j,k\in\Z}$ and $\br,\bn\in\N^2$, $\mathcal
P_\br\,U\,\mathcal P_\bn$ is a non-symmetric finite section of
$U$. We  will write $(\mathcal P_\bn\,U\,\mathcal
P_\bn)^{-1}$ for the  inverse of the symmetric finite section
 $\mathcal P_\bn\,U\,\mathcal P_\bn$ on the finite-dimensional
 subspace  $\mathcal
P_\bn\,\ell^2(\Z)$ (with the understanding that it  cannot  be invertible
on $\ell^2(\Z)$).
\begin{theorem}\label{thm:fs}
Let $\big(\chi (k)\big)_{k\in\Z}$ be a strictly increasing sequence of integers and
$U=(u_{j,k})_{j,k\in\Z}$ be a biinfinite matrix such that  (a) $U^*U$ is
invertible,  and (b)  there exist constants $c,a >0$ such that
\begin{equation}\label{eq:decaying}
\abs{u_{j,k}} \leq c\,e^{-a\abs{j-\chi (k)}}\quad
\text{ for all }j,k\in\Z \, .
\end{equation}
Let $I\subset\N^2$ and assume that for every $\bn\in I$ a finite section $U_\bn := \mathcal P_{\br(\bn)}\,U\,\mathcal P_\bn$
is given such that
\begin{align}\label{uniboundU}
A\,\zweinorm{\boldsymbol{c}}\le \zweinorm{U_\bn\,\boldsymbol{c}}\quad\text{for all }\boldsymbol{c}\in\ell^2(\Z)
\end{align}
for some constant $A>0$ independent of $\bn$.
Then there are constants $\tilde{c},\tilde{a}>0$, such that for all
$\bn\in I$
\begin{equation}
  \label{eq:c1}
\norm{
U\,(U^*\,U)^{-1} \, e_0
- U_\bn\,(U^*_\bn\,U_\bn)^{-1} \, e_0
	}_2 \le \tilde{c}\,e^{-\tilde{a}\,n_0},
\end{equation}
where $n_0:=\min\{n_1,n_2,r_1(\bn),r_2(\bn)\}$.
\end{theorem}
The proof is deferred to the appendix.

Note that the \emph{row} vector  $\big( U (U^*U)\inv ) e_0\big)^* =
e_0^T (U^*U)\inv U^* = e_0^T U^\dagger $ is precisely the zeroth row of the  Moore-Penrose pseudoinverse of
$U$ as it arises in  the computation~\eqref{eq:c15}  of the dual windows $\gamma _L$
and $\gamma ^\circ$.
We also note that the decay condition~\eqref{eq:decaying} models the decay of
the entries off a ``ridge'' $\chi (k)$ rather than off-diagonal decay
(in which case $\chi (k) = k$). The above
definition reflects exactly the
behavior of the pre-Gramian matrix $P_g(x)$ of a window with exponential
decay, since $|P_g(x)_{jk}| = |g(x+\alpha j - k/\beta ) | \leq C e^{-
    a\alpha     |j - k/(\alpha \beta ) |}$, in which case $\chi (k) =
  \lfloor k/(\alpha \beta )\rfloor$. Decay conditions of this type occur in wavelet
  theory~\cite{ABK08} and in the theory of Fourier integral operators~\cite{CGNR13}.

We can now prove the main result of our paper, namely that the
numerically computable dual windows $\gamma _L$ converge exponentially
fast to the canonical dual window $\gamma ^\circ $.

\begin{theorem}\label{thm:conv_gammaL}
The dual windows $\gamma_L$, $L\in\N$, in Algorithm~\ref{algo:gamma}
approximate the canonical dual window $\gamma ^\circ $ of $g$ at an exponential rate
\begin{equation}
  \label{eq:c16}
  \zweinorm{\gamma_L-\gamma^\circ} \le \tilde{c}\,e^{-\rho L}.
\end{equation}
\end{theorem}

\begin{IEEEproof}
  We set  $\bn(L) = (\abs{k_1}+L,k_2+L)$ and $\br(L) =
  (\abs{i_1(L)},i_2(L))$ as in Algorithm~\ref{algo:gamma}. The
  matrices $P_L(x)$ in Algorithm~\ref{algo:gamma} are exactly the
  non-symmetric finite sections  $P_L(x) = \mathcal{P}_{\mathbf{r}(L)} P_g(x)
  \mathcal{P}_{\mathbf{n}(L)}$ of the pre-Gramian $P_g(x)$. Then
  Theorem~\ref{thm:bounds-PL} implies that the  finite sections $P_L(x)$
  of the pre-Gramian $P_g(x)$ are left-invertible (on the appropriate
  finite-dimensional subspaces) with constants independent of
  $L$. Therefore   the decay conditions and the uniform bounds in the
  assumptions of
  Theorem~\ref{thm:fs} are fulfilled. Thus for fixed $x$ we obtain
  that
$$
\norm{e_0^T P_g(x)^\dagger - e_0^T P_L(x)^\dagger }_2 \leq \tilde c
e^{-\tilde a n_0(L)} \, .
$$
Finally the approximation of the dual windows $\gamma _L$ and $\gamma
^\circ$ follows from
\begin{align*}
  \norm{\gamma _L - \gamma ^\circ}_2^2 &= \int _{-\infty } ^\infty
  |\gamma _L (x)- \gamma ^\circ (x) |^2 \, dx  \\
&= \int _0 ^\alpha \sum _{j\in \Z } |\gamma _L (x+\alpha j)- \gamma
^\circ (x+\alpha j) |^2 \, dx \\
& = \beta^2 \int _0^\alpha \norm{e_0^T
  P_L(x)^\dagger - e_0^T P_g(x)^\dagger}_2^2 \, dx  \\
&\leq \alpha(\tilde c \beta)^2  e^{-2\tilde a n_0(L)} \, .
\end{align*}
Since  $n_0(L) :=\min\{n_1(L),n_2(L),r_1(L),r_2(L)\} = \min \{ |k_1|+L, k_2
+ L, |i_1(L)|, i_2(L)\} = L + C$ for some integer constant depending
on the window only, the rate of approximation in \eqref{eq:c16} follows.
\end{IEEEproof}
\begin{remark}
(i)
 Theorem~\ref{thm:fs} is not contained in the results on the
 non-symmetric finite section method in \cite{GroeRzeStroh:2010}. The
 selection of rows by $\br(\bn)$ in our assumption (\ref{uniboundU})
 meets only the condition in \cite[Lemma~5.2]{GroeRzeStroh:2010}, which reads as
$$\sup_{\bn\in I}\norm{(U_\bn^*\,U_\bn)^{-1}}_{\ell^2\rightarrow\ell^2}<\infty$$
in our notation. However, the condition in
\cite[Lemma~5.1]{GroeRzeStroh:2010} is not matched and can only be
satisfied  by considerably increasing the number  of rows of $U_\bn$. The
approximations of the canonical dual $\gamma^\circ$ based on the
non-symmetric finite section method in \cite{GroeRzeStroh:2010} do not
provide dual windows, in contrast to our approximations $\gamma_L$.

(ii)  With an analogous proof, we obtain that the dual windows $\gamma_L$ in Algorithm \ref{algo2:gamma} approximate the canonical dual of the EB-spline $B_\Lambda$ at an exponential rate.
\end{remark}


\section{Discretization and implementation}\label{sec:discrete}

For the numerical use of TP functions and EB-splines in signal
analysis it is often  necessary to discretize these windows. We define the
sampling operator $S_\delta$ for a given sampling rate $\delta>0$ by
$$S_\delta g := \sqrt{\delta}\left( g(\delta k) \right)_{k\in\Z}$$
and the periodization operator $P_K$ with period $K>0$ by
$$P_K g(x) := \sum_{k\in\Z} g(x+Kk),\quad x\in[0,K).$$
For discrete signals $\boldsymbol{c}\in\ell^1(\Z)$ we let
$$P_K\,\boldsymbol{c} := \left(\sum_{k\in\Z}c_{j+Kk}\right)_{j=0,\ldots,K-1}\in\C^K.$$
Furthermore we consider the dilation operator
$$D_h g(x) := \sqrt{h}\, g(hx),\quad h>0,$$
which preserves the $L^2$-norm. The exponential decay
\eqref{TP:decay} of TP functions or the compact support of EB-splines
imply  that 
$S_\delta\, g$ is in $\ell^1(\Z)\subset\ell^2(\Z)$ and $P_Kg\in L^2([0,K))$.
The combination of  both operators yields  finite discrete signals $P_K S_\delta g\in\C^K$.
For $a,b,M,N,K\in\N$ and $Mb = Na =K$, the discrete Gabor system is defined as
$$\begin{array}{l}
\mathcal G(P_K\,S_\delta\,g,a,\tfrac{1}{M}) =
 \left\{ e^{2\pi i\tfrac{l}{M}\cdot}\,P_K\,S_\delta\,g(\cdot-ka)\middle|
 \begin{array}{l}k=0,\ldots,N-1 \text{ and } \\  l=0,\ldots,M-1 \end{array}\right\}.
 \end{array}
$$
The following result is well known and holds for an extremely general
class of window functions.
All
TP functions of finite type and all EB splines
together with the duals of Algorithms \ref{algo:gamma} and \ref{algo2:gamma}
satisfy this mild condition.   

\begin{proposition}\label{prop:Sond}
Let $\alpha,\beta>0$ and
$\alpha\beta =
\tfrac{a}{M} = \tfrac{b}{N}$ and $Mb = Na =K$ with $a,b,M,N,K\in\N$.
Let $g,\gamma\in L^2(\R)$ such that $(g(x+j\alpha))_{j\in\Z}$
and $(\gamma(x+j\alpha))_{j\in\Z}$ are absolutely summable for all $x\in[0,\alpha)$,
and
$$
   \sum_{j\in\Z}\gamma(x+j\alpha)\overline{g(x+j\alpha-k/\beta)} =\beta \delta_{0,k}$$
for all $x\in[0,\alpha)$ and $k\in\Z$.
Then $\mathcal
G(P_K\,S_{\alpha/a}\,g,a,\tfrac{1}{M})$ is a Gabor frame for
$\C^K$ and  $\mathcal G(P_K\,S_{\alpha/a}\,\gamma,a,\tfrac{1}{M})$
is a dual Gabor frame.
\end{proposition}

This statement is proved under slightly different assumptions in \cite{Sond:2007}. The assumptions
in Proposition \ref{prop:Sond} lead directly to the verification of the Wexler-Raz criterion
 in \cite[Theorem A.3]{Sond:2007} for dual Gabor frames of $\C^K$. The details of the proof are omitted here.

\begin{remark}\label{rem:dilation}
Since $S_\delta D_h = S_{\delta h}$, it is helpful to dilate the
function $g$ by the sampling rate. Subsequently we can work with a
sampling rate $\delta = 1$ and consider $P_K S_1 \tilde g =: P_K S
\tilde g$ of some scaled TP function or EB-spline $\tilde
g:=D_{\alpha/a} g$. In many practical situations $a/\alpha$ is
proportional to $\sqrt{K}$. 
Thereby the
time-frequency localization of the window is independent of $K$.
\end{remark}

In the remaining part of this section we describe implementations of
discretized TP functions and EB-splines as well as the duals from
Algorithms~\ref{algo:gamma} and~\ref{algo2:gamma} in Section~\ref{sec:known}. For
this purpose, we use some knowledge about the Zak transform of these
functions. For a parameter $\alpha >0$ and a function $f\in L^2(\R)$
with absolutely summable $(f(x+\alpha j))_{j\in\Z}$ for $x\in\R$,
the Zak transform  is defined by
$$\mathrm{Z}_{\alpha}f(x,\omega) := \sum_{j\in\Z}f(x+\alpha j) e^{-2\pi ij\alpha \omega}.$$
The Zak transform is $\alpha$-quasiperiodic in $x$ and $1/\alpha
$-periodic in $\omega$~\cite{janssen-zak88}. For a given periodization parameter $K\in\N$
we obtain the discrete version of a scaled TP function or EB-spline
$\tilde g$ by
\begin{equation}\label{eq:PSg}
P_K S\tilde g(k) = \mathrm Z_K \tilde g(k,0),\quad k=0,\ldots,K-1.
\end{equation}

\subsection{EB-splines}\label{subsec:EB}
Since EB-splines have compact support, their Zak transform is a finite
sum and only requires finitely many point  evaluations of these functions.

\textsl{Case 1}: The EB-spline $B_{\lambda,\ldots,\lambda}$ with  a
single  weight $\lambda\in\R$ of multiplicity $m\in\N$ can be factorized into an exponential and the cardinal polynomial B-spline $N_m$ of order $m$
\begin{align*}
B_{\lambda,\ldots,\lambda} &= e^{\lambda (\cdot)}\chi_{[0,1)} \, \ast \ldots \ast \, e^{\lambda (\cdot)}\chi_{[0,1)} \\
                              &= e^{\lambda (\cdot)}\left( \chi_{[0,1)} \, \ast \ldots \ast \, \chi_{[0,1)}\right) = e^{\lambda (\cdot)}\, N_m.
\end{align*}
Hence it can be evaluated by the well-known algorithm by Cox and
deBoor (see \cite{deBoor:1978}), which is part of the standard signal
processing toolboxes.

\textsl{Case 2}: For EB-splines with pairwise distinct weights $\lambda_1<\cdots<\lambda_m$ and $m\geq2$ Christensen and Massopust \cite{ChrisMass:2010} give the closed form
\begin{equation}\notag
  B_{\lambda_1,\ldots,\lambda_m}(x+k-1)=\sum_{j=1}^m \alpha_j^{(k)} e^{\lambda_j x},~ x\in[0,1),~1\leq k\leq m,
	\end{equation}
with coefficients
\begin{align*}
\alpha_j^{(k)} = \begin{cases}
\prod\limits_{r=1,\atop r\neq j}^{m}(\lambda_m-\lambda_r)^{-1},\ \ &k=1,\\
 \frac{\sum\limits_{1\leq j_1<\cdots<j_{k-1}\leq m,\atop j_1,\cdots,j_{k-1}\neq j}e^{\lambda_{j_1}+\ldots+\lambda_{j_{k-1}}}}{(-1)^{k-1} \,\prod\limits_{r=1,\atop r\neq j}^{m} (\lambda_m-\lambda_r)},\ \ &k=2,\ldots,m.
\end{cases}
\end{align*}

\textsl{Case 3}: For EB-splines with several distinct weights with multiplicities, we use the following four-term recurrence relation, which was stated in \cite{DynRon:1988} and \cite{Ron:1988}.
\begin{theorem}{\bf \cite{DynRon:1988}, \cite{Ron:1988}}
Let $\lambda_1,\ldots,\lambda_N\in\R$, $\lambda_1\neq\lambda_N$ and
$$B_{\lambda_1,\ldots,\lambda_N} =  e^{\lambda_1(\cdot)}\chi_{[0,1)} \, \ast e^{\lambda_2(\cdot)}\chi_{[0,1)} \,
 \ast \ldots \ast \, e^{\lambda_m(\cdot)}\chi_{[0,1)}$$
the associated EB-spline. Then the following recursion holds
\begin{align*}
 & B_{\lambda_1,\ldots,\lambda_N}(x) = \frac{B_{\lambda_1,\ldots,\lambda_{N-1}}(x)-B_{\lambda_2,\ldots,\lambda_N}(x)}{\lambda_1-\lambda_N}
 +\frac{e^{\lambda_1}B_{\lambda_2,\ldots,\lambda_N}(x-1)-e^{\lambda_N}B_{\lambda_1,\ldots,\lambda_{N-1}}(x-1)}
 {\lambda_1-\lambda_N}.
\end{align*}
\end{theorem}

By iteration of this recurrence it is possible to reduce any given
EB-spline into either several lower order EB-splines with pairwise distinct
weights of multiplicity one or only one weight of higher
multiplicity. These can be treated as in the aforementioned cases.

\subsection{TP functions}
In the case of TP functions, we present two different implementations of the computation of $P_KS \tilde g$ in \eqref{eq:PSg}.
For both  we use that TP functions are invariant under dilation.
\begin{lemma}
Let $g$ be a TP function of finite type $N\in\N$ with weights $(\delta_\nu)_{\nu=1}^N$ and $h>0$ a scaling parameter. Then
$$g_h=\sqrt{h}D_h g $$
 is the TP function of finite type in \eqref{eq:TPFFT} with weights $(\tfrac{\delta_\nu}{h})_{\nu=1}^N$.
\end{lemma}
The first implementation uses the identity in \cite[Remark~2]{BanGroeStoe:2013}
$$\mathrm Z_K{g_h}(x,\omega) = [\tfrac{1}{\delta_1},\ldots,\tfrac{1}{\delta_N}\mid r_{x,\omega}]\qquad \forall x\in[0,K),$$
where the right-hand side is the divided difference of $r_{x,\omega}$ with
$$r_{x,\omega}(y) = (-1)^{N-1}h\,\left(\prod_{\nu=1}^N \delta_\nu^{-1}\right)\,
\frac{e^{-hxy}}{1-e^{-K(hy+2\pi i\omega)}}$$
in the knots $\delta_1^{-1},\ldots,\delta_N^{-1}$.

The second  implementation uses the connection of TP functions to
EB-splines. In \cite[Theorem~3.4]{KloStoe:2014} it is shown  that the Zak transform of a TP function can be expressed in terms of the Zak transform of an associated EB-spline.
\begin{theorem}{\bf \cite{KloStoe:2014}}\label{b-splinezak}
Let $g$ be a TP function of finite type with weights $\delta_1,\ldots,\delta_N\in\R$.
With $\lambda_{\nu}:=-\frac{Kh}{\delta_{\nu}}$, $\nu=1,\ldots,N$, we have
$$K\mathrm Z_K g_h(x,0) = \prod_{\nu=1}^N \, \frac{\lambda_\nu}{e^{\lambda_\nu}-1}\,
\mathrm Z_1B_{\lambda_1,\ldots,\lambda_N}(\tfrac{x}{K},0),\ \ x\in[0,1).$$
\end{theorem}
Consequently $P_KS \tilde g$  can be computed by the Zak transform of the corresponding EB-spline as described in \ref{subsec:EB}.

\subsection{Dual windows}\label{imp:duals}
For TP functions $g$ and lattice parameters $\alpha,\beta>0$ with $\alpha\beta<1$, Algorithm~\ref{algo:gamma} in Section~\ref{sec:known} allows us to compute samples $S_{\alpha/N}\gamma_L$, $N\in\N$, of a dual window $\gamma_L$. Likewise, Algorithm \ref{algo2:gamma} provides the sampled dual windows of EB-splines $B_\Lambda$.
Therefore, for a given periodization parameter $K\in\N$ and the time-shift parameter $a\in\N$, we compute
$$P_KS_{\alpha/a}\gamma_L(k)= \sum_{j\in\Z}
(S_{\alpha/a}\gamma_L)(k+jK),\quad k=0,\ldots,K-1 \, .$$
This sum is  finite because of  the compact support of $\gamma_L$.

\begin{center}
\begin{figure}
\hspace*{4.5cm}%
\includegraphics[width=8.8cm]{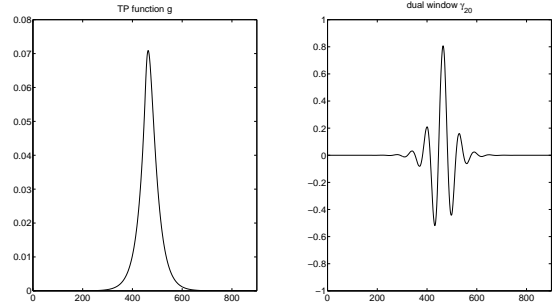}

\caption{Discrete TP function and its dual $\gamma_{20}$ for $a=20$ and $M=30$.}\label{plot:discretedual}
\end{figure}
\end{center}

\begin{example}
We use Algorithm~\ref{algo:gamma} for the computation of the dual Gabor window $\gamma_L$ of the asymmetric TP function $g$ with parameters $\delta=[-1,1,1/3,1/5]$ and lattice parameters $\alpha=2/3$, $\beta=1$. The discretization parameters $K=900$, $a=20$, $b=1/30$, are chosen according to the standard dilation $a/\alpha=\sqrt{K}$ in Remark \ref{rem:dilation}. Figure~\ref{plot:discretedual} shows the discrete TP function $P_K S_{\alpha/a}g$ and its dual window $P_K S_{\alpha/a} \gamma_L$ for $L=20$. The difference $\norm{P_K S_{\alpha/a} (\gamma_L-\gamma^\circ)}_2$ to the discrete canonical dual is $7\cdot 10^{-8}$ measured in the $\ell_2$-norm of $\C^{900}$.
\end{example}

Sometimes  the frame bounds in the discrete case may  be better than
in the continuous case. Therefore the discrete TP functions or
EB-splines may  even provide Riesz bases at the critical density, as is
explained in the  following recent
result.

\begin{proposition}{\bf \cite{BanGroeStoe:2013},\cite{Klo:2015}}
Let $g$ be a continuous TP function as in \eqref{TPall}, including the infinite type with $\eta = 0$. Assume $\alpha = M\in\N$ and let $\beta = 1/M$ and $K\in\N$ such that $K/M\in\N$.
If $K/M$ is odd, then $\mathcal G(P_KSg,\alpha,\beta)$ is a basis of $\mathbb C^K$.

In addition, assume that $g$ is even, which means that $\{\delta_{\nu}\mid \delta_{\nu}>0\}=\{-\delta_{\nu}\mid \delta_{\nu}<0\}$.
If $M$ is odd, then $\mathcal G(P_KSg,\alpha,\beta)$ is a basis of $\mathbb C^K$.
\end{proposition}

\begin{remark}
In the critical case $\alpha \beta =1$ the computation of the dual window $\gamma$
cannot be performed by the aforementioned discretization procedure, as
the Gabor system $\gab$ is not a frame in $L^2(\R)$. In this case,
 the usual method using the discrete Fourier transform and the
discrete Zak transform of $P_K S\tilde g$ should be applied for the
computation of the discrete dual window~\cite{BoelHlaw97}.
\end{remark}

\section*{Acknowledgement}
First discussions on the topic of this article were performed while J. St\"ockler visited the Erwin Schr\"odinger Institute. The authors are grateful to Z.~{P}r\r{u}\v{s}a and P.~L. {S}{\o}ndergaard for many helpful discussions concerning the implementation of the algorithms for the computation of the canonical dual window which were held at the Strobl conference on Modern Time-Frequency Analysis in 2014.

\section*{Appendix}
In the appendix we provide the technical details of the proofs of the
main theorems. The proof of Theorem~\ref{thm:bounds-PL} depends heavily on
\cite[Theorem~9]{GroeStoe:2013}.

\begin{IEEEproof}[Proof of Theorem~\ref{thm:bounds-PL}]
The pre-Gramian is bounded on $\ell ^2(\Z )$ as a consequence of the
exponential decay of $g$ and Schur's test. Therefore the finite
sections $P_L(x) = \mathcal P_{\br(\bn)}\,P_g(x) \,\mathcal P_\bn$ are
uniformly bounded.

To prove the existence of a lower bound $A$, we construct a left
inverse of $P_L(x)$  by adapting the proof of
\cite[Theorem~9]{GroeStoe:2013}.
By step~2 of Algorithm~\ref{algo:gamma}  $P_L(x)$ has columns indexed by $k,
k_1-L\le k\le k_2+L$ and rows indexed by $j, i_1(L)\le j\le i_2(L)$,
thus every  left-inverse  $Q_L(x)$ has columns indexed by $i_1(L)\le
j\le i_2(L)$ and rows indexed by $k_1-L\le k\le k_2+L$. We will
construct these  rows one by one in three steps.
\medskip

\noindent {\emph Step 1.} For every index $-L\le \ell \le L$, we choose the following submatrix $R_\ell$
of $P_L(x)$. Given $x$ and $\ell \in \mathbb{Z}$, we choose  $y_\ell\in [0,\alpha)$ and $j_\ell\in\Z$, such that
\[
    y_\ell = x+j_\ell\alpha-\frac{\ell}{\beta}.
\]
Then we define the matrix $R_\ell=P_0(y_\ell)$ by using the
same definitions as in step~3 of Algorithm \ref{algo:gamma} with $L=0$;
i.e., we let
\begin{align*}
     i_1^\ell(0) &= \left\lfloor \frac{k_1 + m -
        1}{\alpha\beta}-\frac{y_\ell}{\alpha}\right\rfloor+1,\\
     i_2^\ell(0) &=
    \left\lceil \frac{k_2 - n +
        1}{\alpha\beta}-\frac{y_\ell}{\alpha}\right\rceil-1,
 \end{align*}
and
\[
   R_\ell=\left( r_{\ell;j,k}\right)
	_{i_1^\ell(0)\le j\le i_2^\ell(0),~k_1\le k\le k_2},
\]
with $r_{\ell;j,k}=g\left(y_\ell+j\alpha-\frac{k}{\beta}\right)$. Note that
\[
   y_\ell+j\alpha-\frac{k}{\beta} = x+(j+j_\ell)\alpha-\frac{\ell + k}{\beta}.
\]
Hence,   the columns of $R_\ell$, as indexed by $k_1\le k\le k_2$, correspond to sections of
 columns of $P_L(x)$ indexed by $\ell+k_1\le k\le \ell+k_2$.
 More precisely, $R_{-L}$ contains sections of the first
$k_2-k_1+1$ columns of $P_L(x)$;  with increasing $\ell$ the selection
of columns is shifted to the right;  and $R_L$ contains sections of the
last $k_2-k_1+1$ columns of $P_L(x)$. For the rows of $R_\ell$, we
observe that
\[
   i_1(L)\le \left\lfloor \frac{k_1 + m -
        1}{\alpha\beta}-\frac{x}{\alpha}+ \frac{\ell}{\alpha\beta}\right\rfloor+1
        =i_1^\ell(0)+j_\ell
\]
and $i_2(L)\ge i_2^\ell(0)+j_\ell$. To summarize, each $R_\ell $  is a $(i_2^\ell (0) -
i_1^\ell (0) +1) \times (k_2 - k_1+1)$-submatrix of $P_L(x)$ with a
dimension independent of  $|\ell | \leq L$. 
\medskip

\noindent{\emph Step 2.} The arguments in \cite[Theorem~8]{GroeStoe:2013}
relate left-inverses of $R_\ell=P_0(y_\ell)$ to left-inverses
of the biinfinite matrix $P_g(y_\ell)$. It is shown that
\begin{itemize}
\item $R_\ell$ has full column rank,
\item there exists a uniform bound $C>0$, which does not
depend on  $y_\ell\in[0,\alpha)$, and left-inverses
$\Gamma_\ell$ of $R_\ell$ such that
\[
   \| \Gamma_\ell\| \le C\quad\mbox{for all}\quad
	-L\le \ell\le L,
\]
\item the rows with index $k_1\le k\le 0$ of $\Gamma_\ell$
are orthogonal to all columns with index $k'>k_2-n$ of the biinfinite matrix $P_g(y_\ell)$, and hence orthogonal to all columns $k'>\ell+k_2-n$ of $P_L(x)$.
Likewise, rows $0\le k\le k_2$ of $\Gamma_\ell$
are orthogonal to all columns with index $k'<\ell+k_1+m$ of $P_L(x)$.
\end{itemize}

\noindent{\emph Step 3.}
With these properties, we obtain the left-inverse $Q_L(x)$ of $P_L(x)$ by
defining the rows of $Q_L(x)$ as follows:
\begin{itemize}
\item We start with $\Gamma_{-L}$ and take its rows $k_1\le k\le 0$ as the first
 rows of $Q_L(x)$, extended by zeros such that
\begin{align*}
     &(Q_L(x))_{k-L,j+j_{-L}}
     = \begin{cases} (\Gamma_{-L})_{k,j},&~  i_1^{-L}(0)\le
		j\le i_2^{-L}(0),\\
		0,&~\mbox{otherwise}.
		\end{cases}
\end{align*}
\item For $-L+1\le \ell\le L-1$ we take the row
with index $0$ of $\Gamma_\ell$ and extend this row by
zeroes,
\[
     (Q_L(x))_{\ell,j+j_\ell}= \begin{cases} (\Gamma_{\ell})_{0,j},&~ i_1^\ell(0)\le
		j\le i_2^\ell(0),\\
		0,&~\mbox{otherwise}.
		\end{cases}
\]
\item We end with $\Gamma_{L}$ and take its rows $0\le k\le k_2$ as the last
 rows of $Q_L(x)$, extended by zeros such that
\begin{align*}
     &(Q_L(x))_{k+L,j+j_{L}}
     = \begin{cases} (\Gamma_{L})_{k,j},&~
     i_1^L(0)\le j\le i_2^L(0),\\
		0,&~\mbox{otherwise}.
		\end{cases}
\end{align*}
\end{itemize}
It is clear that all entries of $Q_L(x)$ are bounded  by
the constant $C$ in step 2. Moreover,
every row and column of $Q_L(x)$ has at most
\[
  \max_{-L\le \ell\le L}  (i_2^\ell(0)-i_1^\ell(0)+1 ) \le \frac{k_2-k_1+1}{\alpha\beta}
\]
nonzero entries. Therefore, by Schur's test, we obtain that
\[
  \|Q_L(x)\| \le C\,\frac{k_2-k_1+1}{\alpha\beta}.
\]
Since $Q_L $ is a left inverse of $P_L$, we obtain
$$
\norm{c}_2 = \norm{Q_L(x) P_L(x) c}_2 \leq C \,
\frac{k_2-k_1+1}{\alpha\beta}  \norm{P_L(x)c}_2 ,
$$
and we may choose $A = (C \,
\frac{k_2-k_1+1}{\alpha\beta} )\inv $ as a lower bound in~\eqref{PLbounds}.
This completes the proof of (\ref{PLbounds}). We conclude that
$\zweinorm{e_0^T\,P_L(x)^\dagger}\leq C \,
\frac{k_2-k_1+1}{\alpha\beta} =  A^{-1}$. Therefore
\begin{align*}
  \zweinorm{\gamma_L}^2 &= \int _{-\infty } ^\infty |\gamma _L(x)|^2
  \, dx \\
& = \beta^2\, \int_0^\alpha \sum _{j\in \Z } |\gamma _L(x+j\alpha )|^2 \,dx \\
&= \beta ^2 \int _0^\alpha  \zweinorm{e_0^T\,P_L(x)^\dagger}^2\, dx \leq \alpha\,\beta^2\,A^{-2}.
\end{align*}
\end{IEEEproof}

In the following we will need a well known result of
Jaffard~\cite{Jaff:1990} and Baskakov~\cite{Bas90}.

\begin{proposition} \label{jaffexp}
  Let $A = (A_{jk})_{j,k \in \Z } $ be a biinfinite matrix with
  exponential off-diagonal decay, i.e., there exist  constants $C,a>0$, such that
  \begin{equation}
    \label{eq:c89}
|A_{jk}| \leq C e^{-a |j-k|} \quad \forall j,k \in \Z \, .
  \end{equation}
If $A$ is invertible, then $A\inv $ also possesses exponential
off-diagonal  decay
and there exist $C'> 0$ and $0< \tilde a < a$, such that
$$
|(A\inv )_{jk}| \leq C' e^{-\tilde a |j-k|} \qquad \forall j,k\in \Z
\, .
$$
\end{proposition}

\begin{corollary} \label{jafffsm}
  Let $A$ be  positive and invertible on $\ell ^2(\Z )$ with
  exponential off-diagonal decay~\eqref{eq:c89} and let $\mathcal{P}_{\bn }A \mathcal{P}_{\bn } $ be a
  sequence of finite sections. Then every matrix  $\mathcal{P}_{\bn }A
  \mathcal{P}_{\bn }$ is invertible on $\mathcal{P}_{\bn}\ell ^2(\Z)$
  and there exist $C'>0$ and $\tilde a,
  0<\tilde a <a$, such that
  \begin{equation}
    \label{eq:c17}
|(\mathcal{P}_{\bn }A  \mathcal{P}_{\bn })\inv_{jk}|\leq C' e^{-\tilde a
  |j-k|} \qquad  \text{ for } -n_1 \leq j,k \leq n_2 \, .
  \end{equation}
\end{corollary}
\begin{IEEEproof}
  The proof follows~\cite{HagRochSilb:2001} and \cite{GroeRzeStroh:2010}. Let $\mathcal{A}$ be the matrix obtained
  by stacking the finite sections $\mathcal{P}_{\bn }A
  \mathcal{P}_{\bn }$ along the diagonal. Then $\mathcal{A}$ possesses
  exponential off-diagonal decay~\eqref{eq:c89}. Since $A$ is invertible and positive, its  spectrum
  $\sigma (A)$ is
  contained in an interval $[A,B]$ for some $A,B>0$, consequently  the spectrum of the
  restriction of $\mathcal P_\bn\,A\,\mathcal P_\bn$ on
  $\mathcal P_\bn\ell^2(\Z)$ is also contained in $[A,B]$ and every
  finite section is invertible  on $\mathcal P_\bn \ell ^2(\Z
  )$. Therefore the stacked matrix $\mathcal{A}$ is invertible on
  $ \oplus \mathcal P_\bn\ell^2(\Z) \simeq \ell ^2(\Z )   $. By
  Proposition~\ref{jaffexp} $\mathcal{A}\inv $ possesses exponential
off-diagonal   decay. Since $\mathcal{A}\inv $ consists of the blocks
  $(\mathcal{P}_{\bn }A  \mathcal{P}_{\bn })\inv$, they satisfy
$|(\mathcal{P}_{\bn }A  \mathcal{P}_{\bn })\inv_{jk}|\leq C' e^{-\tilde a
  |j-k|}$ for $ -n_1 \leq j,k \leq n_2 $.
\end{IEEEproof}
We remark that clearly every finite matrix possesses exponential
off-diagonal decay. The point is that the constants may be chosen
independently of the size of the finite section.

\begin{IEEEproof}[Proof of Theorem~\ref{thm:fs}]
Recall that, for   $\boldsymbol{n} =
(n_1,n_2) \in\N^2$, $\mathcal P_{\boldsymbol{n}}\,b =
(\ldots,0,b_{-n_1},b_{-n_1+1},\ldots,b_{n_2-1},b_{n_2},0,\ldots)^T$   is the orthogonal projection onto
$\mathcal P_{\boldsymbol{n}}\,\ell^2(\Z) \cong \C^{n_1+n_2+1}$. We
write  $U_\bn= \mathcal
P_\br\,U\,\mathcal P_\bn$ for a non-symmetric finite section of
$U$.  In the assumption of Theorem~\ref{thm:fs} $\br $ depends on
$\bn$, but we will omit this dependence in the notation. All operator norms are in $\ell^2(\Z)$.

Let $b:= U\,(U^*\,U)^{-1} \, e_0$ and $d_\bn := U_\bn\,(U^*_\bn\,U_\bn)^{-1}
\, e_0$. We decompose the norm into three parts
\begin{align}\notag
\zweinorm{b - d_\bn} &\le \zweinorm{(U - U\,\mathcal P_\bn)\,(U^*\,U)^{-1}\,e_0}\\ \label{pf:finsec}
+& \zweinorm{U\,\mathcal P_\bn\,\Big((U^*\,U)^{-1} - (\mathcal P_\bn\,U^*\,U\,\mathcal P_\bn)^{-1}\Big)\,e_0}\\ \notag
+& \zweinorm{U\,\mathcal P_\bn\,((\mathcal P_\bn\,U^*\,U\,\mathcal P_\bn)^{-1} - U_\bn\,(U^*_\bn\,U_\bn)^{-1})\, e_0}.
\end{align}
\medskip

\noindent{\emph 1.} Note that $U^*\,U$ is invertible and positive. Since $U$
  fulfills the decay property (\ref{eq:decaying}), it is easy to see
  that  the symmetric
  matrix $U^*\,U$ decays exponentially off the diagonal, i.e., for
  some constants $C,a>0$
$$
|(U^*U)_{jk}| \leq C e^{-a |j-k|} \qquad j,k \in \Z \, .
$$
By Proposition~\ref{jaffexp} 
the inverse matrix inherits the exponential decay, and thus  there exist  constants $C'
>0$  and   $\tilde a , 0<  \tilde a < a$, 
  such that $|(U^*\,U)^{-1}_{jk}|  \leq   C' e^{-\tilde a |j-k|} $ for all $j,k \in \Z $.
Hence the entries of the vector $v:=(U^*\,U)^{-1}\,e_0$ also decay
exponentially as 
\begin{equation}\label{vj_decay}
\abs{v_j} \le c_1\,e^{-\tilde{a}\,\abs{j}}  \quad j\in \Z \,. 
\end{equation}
Since
\begin{align*}
\zweinorm{(I-\mathcal{P}_\bn ) v}^2 &= \sum _{j= -\infty } ^{-n_1-1}
|v_j|^2 + \sum _{j= n_2+1 } ^{\infty} |v_j|^2
\leq C' \sum _{|j|> n_0} e^{-2\tilde a |j|} = \mathcal{O}(e^{-2\tilde a
  |j|})\, ,
\end{align*}
the decay property \eqref{vj_decay} implies that
\begin{align*}
\zweinorm{(U - U\,\mathcal P_\bn)\,(U^*\,U)^{-1}\,e_0} &\le
\norm{U}\,\zweinorm{(I - \mathcal P_\bn)\,v}
\le
c_2\,e^{-\tilde{a}\,n_0}.
\end{align*}
\medskip

\noindent{\emph 2.}  Since $U^*\, U$ is invertible and positive with   spectrum
  $\sigma (U^*U) \subseteq [A,B]$ for $A,B>0$,   the spectrum of the
  finite sections  $\mathcal P_\bn\,U^*\,U\,\mathcal P_\bn$ on
  ${\mathcal P_\bn\ell^2(\Z)}$ is also contained in $[A,B]$.  As in the finite section method in \cite{HagRochSilb:2001} and with $v=(U^*\,U)^{-1}\,e_0$, we obtain
\begin{align*}
&\zweinorm{((U^*\,U)^{-1} - (\mathcal P_\bn\,U^*\,U\,\mathcal P_\bn)^{-1})\,e_0}
=\zweinorm{(\mathcal P_\bn\,U^*\,U\,\mathcal P_\bn)^{-1}\,\mathcal P_\bn\,U^*\,U\,(\mathcal P_\bn-I)\,v}
\le \norm{(U^*\,U)^{-1}}\,\norm{U^*\,U}\,\zweinorm{(\mathcal P_\bn - I)\,v}
\end{align*}
and consequently
\begin{align*}
\zweinorm{U\,\mathcal P_\bn\,\big((U^*\,U)^{-1} - (\mathcal
  P_\bn\,U^*\,U\,\mathcal P_\bn)^{-1}\big)\,e_0}&
  \le c_3\,e^{-\tilde{a}\,n_0}&.
\end{align*}
\medskip

\noindent{\emph 3.} To treat the third term in~\eqref{pf:finsec},   we need a
  geometric interpretation of the rows of the Moore-Penrose
  pseudoinverse. Let
$$b_\bn := (U\,\mathcal P_\bn)\,(\mathcal P_\bn\,U^*\,U\,\mathcal P_\bn)^{-1}\,e_0$$
and $d_\bn$ as above. Then $b_\bn$ is  the transpose of the zeroth row
of the Moore-Penrose pseudoinverse of $U\,\mathcal P_\bn$. By
Corollary~\ref{jafffsm} $(\mathcal P_\bn\,U^*\,U\,\mathcal
P_\bn)^{-1}$ satisfies \eqref{eq:c17} independently of $\bn$. Therefore
the same argument as in the first part implies  the decay property
\begin{align}\label{dec:b_n}
\abs{(b_\bn)_j} \le c_4\,e^{-\tilde{a}\,\abs{j}}.
\end{align}
It is
essential  that the constants are independent of $\bn$.
Since the Moore-Penrose pseudoinverse of $U\,\mathcal P_\bn$ is also a
left-inverse, we have
$$
b_\bn\,\bot\,V_\bn := \mathrm{span}\{u_k\mid -n_1\le k\le n_2,\ k\neq 0\} ,
$$
where $u_k$, $k\in\Z$, are the columns of the matrix $U$.
Likewise for $U_\bn = \mathcal{P}_{\br} U \mathcal{P}_\bn $ we have
\begin{equation*}
d_\bn\,\bot\,\mathcal P_\br V_\bn\, .
\end{equation*}

Using this orthogonality, we rewrite the vectors $b_\bn $ and $d_\bn $
as follows. Let $\Pi_W$ denote  the orthogonal projection onto some
subspace $W$. Now  set
\begin{equation*}
  \tilde{b}_\bn := (I - \Pi_{V_\bn})\, u_0\quad \text{ and }\quad
  \tilde{d}_\bn := (I - \Pi_{\mathcal P_\br V_\bn})\, \mathcal
  P_\br\,u_0 \, .
\end{equation*}
Since $b_\bn \in \mathrm{Im} (U \mathcal{P}_\bn )$ and $d_\bn \in
\mathrm{Im} (U_\bn )$, we obtain
\begin{equation}\label{normal}
b_\bn = \frac{\tilde{b}_\bn}{\norm{\tilde{b}_\bn}_2^2} \quad \text{
  and } \quad d_\bn =
\frac{\tilde{d}_\bn}{\norm{\tilde{d}_\bn}_2^2}\, .
\end{equation}
 The normalization in (\ref{normal}) is obtained from
\begin{align*}
&\sprod{b_\bn}{u_0} = 1 \quad\text{and}
\quad \sprod{\tilde{b}_\bn}{u_0} = \sprod{\tilde{b}_\bn}{\tilde{b}_\bn + \Pi_{V_\bn}\,u_0} = \norm{\tilde{b}_\bn}_2^2,\\
&\sprod{d_\bn}{\mathcal P_\br u_0} = 1 \quad\text{and}
\quad \sprod{\tilde{d}_\bn}{\mathcal P_\br u_0} = \sprod{\tilde{d}_\bn}{\tilde{d}_\bn + \Pi_{\mathcal P_\br V_\bn}\,\mathcal P_\br u_0} = \norm{\tilde{d}_\bn}_2^2.
\end{align*}
As the third term   in (\ref{pf:finsec}) equals $\norm{b_\bn-d_\bn}_2$, we
first consider $\norm{\tilde{b}_\bn - \tilde{d}_\bn}_2$. For this
purpose we write
\begin{align} \label{hm1}
\tilde{b}_\bn - \tilde{d}_\bn &= (I-\mathcal P_\br)\,(u_0
-\Pi_{V_\bn}\,u_0) + \Pi_{\mathcal P_\br V_\bn}\,\mathcal P_\br\,u_0
- \mathcal P_\br\,\Pi_{V_\bn}\,u_0 \notag \\
                              &= (I-\mathcal P_\br)\,\tilde{b}_\bn +
                              \Pi_{\mathcal P_\br V_\bn}\,\mathcal
                              P_\br\,(u_0 - \Pi_{V_\bn}\,u_0) \notag\\
                              &=    (I-\mathcal P_\br)\,\tilde{b}_\bn
                              +  \Pi_{\mathcal P_\br
                                V_\bn}\,\mathcal P_\br\,\tilde{b}_\bn.
\end{align}
By the assumption (\ref{uniboundU}), the truncated columns $\mathcal
P_\br\,u_k$, with $-n_1\le k\le n_2$ and $k\neq 0$, form a Riesz basis
for $\mathcal P_\br V_\bn$ with lower Riesz bound $A$. Furthermore
it holds that
\begin{align*}
\zweinorm{\Pi_{\mathcal P_\br V_\bn}\,\mathcal P_\br\,\tilde{b}_\bn}^2
&\le A^{-2}\, \sum_{-n_1\le k\le n_2,\,k\neq 0}\abs{\sprod{\Pi_{\mathcal P_\br V_\bn}\,\mathcal P_\br\,\tilde{b}_\bn}{\mathcal P_\br\,u_k}}^2\\
&= A^{-2}\, \sum_{-n_1\le k\le n_2,\,k\neq 0}\abs{\sprod{\mathcal P_\br\,\tilde{b}_\bn}{u_k}}^2\\
&= A^{-2}\, \sum_{-n_1\le k\le n_2,\,k\neq 0}\abs{\sprod{(\mathcal P_\br-I)\,\tilde{b}_\bn}{u_k}}^2\\
&\le A^{-2}B\, \zweinorm{(\mathcal P_\br-I)\,\tilde{b}_\bn}^2,
\end{align*}
where $B=\norm{U}^2$ denotes the Bessel bound of all columns $u_k$,
$k\in\Z$.
Taking the $\ell ^2$-norm in \eqref{hm1} and substituting the above
estimate, we obtain
$$\norm{\tilde{b}_\bn - \tilde{d}_\bn}_2 \le c_5\,\norm{(I-\mathcal P_\br)\,\tilde{b}_\bn}_2.$$
We now return to $\norm{b_\bn-d_\bn}_2$. It is an easy exercise that
for nonzero vectors $y,w$ we have
$$\zweinorm{\frac{y}{\zweinorm{y}^2} - \frac{w}{\zweinorm{w}^2}} \le \frac{3\,\zweinorm{y-w}}{\min\{\zweinorm{y}^2,\zweinorm{w}^2\}}.$$
Using once more that $b_\bn \in \mathrm{Im} (U \mathcal{P}_\bn )$ and $d_\bn \in
\mathrm{Im} (U_\bn )$,  the assumption (\ref{uniboundU}) implies
that
$$0<A\le\norm{\tilde{b}_\bn}_2,
\norm{\tilde{d}_\bn}_2\le\zweinorm{u_0} \, .$$
 By the decay property (\ref{dec:b_n}) we obtain
\begin{align*}
\norm{b_\bn - d_\bn}_2\le\tfrac{3}{A^2}\, \norm{\tilde{b}_\bn - \tilde{d}_\bn}_2&\le\tfrac{3\,c_5}{A^2}\, \norm{(I-\mathcal P_\br)\,\tilde{b}_\bn}_2
\le c_6\,e^{-\tilde{a}\,n_0}.
\end{align*}
To finish the proof of Theorem~\ref{thm:fs},  we add the three contributions in \eqref{pf:finsec} and obtain the
convergence rate \eqref{eq:c1} with a constant  $\tilde{c} = c_2+c_3+c_6$.
\end{IEEEproof}

\end{document}